\documentclass{article}

\usepackage[utf8]{inputenc}
\usepackage{authblk}
\usepackage{xspace}
\usepackage[comma]{natbib}
\usepackage{epigraph} 
\usepackage{graphicx}
\usepackage{booktabs}
\usepackage{makecell}
\usepackage{amssymb}
\usepackage{amsmath}
\usepackage{enumitem}
\usepackage{url}

\newcommand{\BU}{Bucknell University\xspace}
\newcommand{\CHE}{Inside Higher Ed\xspace}
\newcommand{\GIT}{\url{https://github.com/tm032/final_exam_scheduler/tree/main}~\citep{github}\xspace}
\newcommand{\beginFEappendix}{\appendix}
\newcommand{\finishFEappendix}{ }
\newcommand{\tableFont}{\footnotesize}

\title{Final Exam Scheduling at \BU: \\ A Case Study and Open-Source Tool}

\author[1]{Clara Chaplin}
\author[1]{Stanley Gai}
\author[1]{Samuel C. Gutekunst}
\author[1]{Tsugunobu Miyake}
\author[2]{Thiago Serra}
\author[1]{Luke Snyder}
\author[1]{Vy Tran}
\author[1]{Lucas Waddell}

\affil[1]{Bucknell University}
\affil[2]{University of Iowa}

\date{September 2025}

\begin{document}

\maketitle

\begin{abstract}
\noindent
\textbf{Problem Definition:} Final exam scheduling is a common but challenging optimization problem.   At \BU, a small liberal arts institution, the problem is particularly complex and has historically required the Registrar's Office to spend months manually designing an exam schedule each semester.

\noindent
\textbf{Methodology:} We worked in close collaboration with the Registrar's Office. First, we created visualization tools to help their manual scheduling process.  Then we designed integer programming models and heuristics to produce a portfolio of possible exam schedules.  Finally, we developed open-source, user-friendly software, enabling their office to directly produce and adjust these schedules.

\noindent
\textbf{Results and Managerial Implications:} Our tools -- both software and algorithms -- are now in use at \BU.  This collaboration has led to both substantial time savings and improved schedules.  Since the implementation of this project, for example, the proportion of students who have a back-to-back exam in a given semester has decreased from roughly a third to about 10\%.  Our tools are fully open-source and rely on an open-source optimization solver, and our approach garnered national media attention with an article in \CHE.

\end{abstract}

\section{Introduction}\label{sec:Intro}

Scheduling problems are ubiquitous in higher education, from planning course schedules~\citep{Gon18,Chr24} to scheduling medical internships~\citep{Akb22}.  Perhaps the most notorious problem, however, is final exam scheduling. 
The literature on this topic includes, among others, 
case studies for Cedar Crest College~\citep{Meh81}, 
Universidad Polit{\'e}cnica de Madrid~\citep{Rom82}, 
Technion~\citep{Str17}, 
secondary education in Norway~\citep{Ave22}, 
and Cornell University~\citep{Ye24}, 
as well as an examination track at the International Timetabling Competition~\citep{Mcc07,Mcc12}.

One canonical approach is to assign exam times based on when classes meet during the semester. 
That prevents most, but not necessarily all, students from having to take two exams at the same time: Unlike Hermione Granger in the Harry Potter series~\citep{JKR99}, our students are not able to travel in time to attend simultaneous lectures. However, they may sometimes be granted an exception and allowed to attend two courses whose meeting times overlap. In addition, that approach to exam scheduling may not be tractable for schools which offer multiple sections of the same course at different meeting times, such as for calculus and introductory physics sequences.  The final exam time for all such sections may need to be coordinated to ensure that all students are assessed with the same final exam at the same time.  Having separate exam times for every class meeting time and for each of these larger multi-section courses can be difficult.  For instance, to give each class time and each class requesting a coordinated exam slot at \BU its own exam slot, the final exam period would need to be three times longer.

Exams can become more difficult to schedule at liberal arts institutions, where sections of courses are intentionally kept small; there may be proportionately more of these sections offered at different times and the number of large classrooms on campus is limited. 
For example, at \BU there are dozens of sections of calculus each semester with different meeting times.  Another challenge in liberal arts institutions is that university-level requirements encourage students to take courses far from their home departments.  While it may be possible for each college, department, or major of some universities to independently create a final exam schedule, that approach is not tenable at liberal arts institutions due to frequent cross-department enrollment.

On top of avoiding overlapping exams, we also want to produce exam schedules that support student success.   Final exams tend to be several hours. 
If a student has a short or nonexistent interval between exams, 
that affects their ability to properly sleep, eat,  and prepare for their next exam.

Finally, every institution has their unique rules and constraints.  For example, \BU  has multiple colleges, and the final exam schedule must be approved by the deans of each college before the start of the academic year. 
Consequently,  final exams are scheduled before the add-and-drop period during which course enrollments are finalized.
Similarly, as an an NCAA D1 College, 
some departments require that exams avoid  days when student-athletes might have to travel.

In this paper, we focus on exam scheduling at \BU where exam scheduling is handled by the Registrar's Office, and where previously staff members would  manually build an exam schedule over multiple months. (Hereafter, we use \emph{registrar} as shorthand for the staff in the Registrar's Office who build the exam schedule.)  We developed integer programming models for exam scheduling and created open-source software allowing the registrar to create and modify optimized schedules.  In doing so, our main contributions are threefold.

First, our work has led to substantial improvements at \BU.  We have greatly reduced the number of students who have an \emph{exam inconvenience}: two exams at the same time, three exams within 24 hours, back-to-back exams, a night exam followed by a morning exam the next day, or four exams within 48 hours.  With just under 4,000 students, prior to our work 1,182 students --- close to 30\% --- had at least one exam inconvenience in Spring 2023. After our optimization models were used, only 432 students had an inconvenience in Spring 2024. 

Second, our work provides a roadmap for applications of Operations Research and Analytics in academic settings.  We worked in close collaboration with the Registrar's Office, first learning about their manual process and building visualization tools to support it.  Concurrently, we surveyed stakeholders on what aspects of exam scheduling mattered most to them, including a formal survey of hundreds of students.  Then, we introduced integer-programming based models to provide the registrar with a portfolio of improved schedules.  Finally, we built fully open-source software tailored to their process: our software links with automatically-updating enrollment data; allows the registrar to ``hard-code'' the requests that they accept from instructional staff, like requiring a specific exam to avoid certain days; generates a portfolio of schedules with a click of a button; and then allows the registrar to ``drag and drop'' to make any final adjustments, seeing live-updates on how those changes affect the schedule.  The registrar uses this software through a web-browser interface, and no computer programming knowledge is required on their end.

Finally, we contribute to the broad literature on exam scheduling and bridge what we perceive as a divide within it.  
On the one hand, 
much of the literature focuses on  highly-tailored models for specific institutions.  
On the other hand, 
another major strand has been methodological and focused on stylized and general exam scheduling problems. 
While our algorithms and tools are motivated by the specifics of exam-scheduling at \BU,  we designed these tools to be as generalizable as possible for other institutions: they are open-source, highly flexible, easily tailorable, and fully rely on open-source optimization solvers.  

The rest of this paper is organized as follows. 
We provide a brief literature review in Section~\ref{sec:litrev}.  Then, in Section~\ref{sec:atBU}, we formalize the problem at \BU and highlight our survey data on student exam preferences.  We provide more context on our approach, including our visualization tools, integer programming model, and software in Section~\ref{sec:Approach}. We conclude by highlighting the results of this project at \BU in Section~\ref{sec:Results}.

\section{Related Work: Good Luck Anyway}\label{sec:litrev}

\epigraph{\textit{Having tried to do in the past what you are proposing now, I hope you have the time and resources to match your enthusiasm.  Do \underline{not} underestimate the complexity of the problem!  Good luck anyway.}}{Comment from a respondent of a survey of Registrars conducted by \cite{Bur96b}.}

There is a vast literature on exam scheduling dating back to at least the early 1960's, when  \cite{Bro64} provided a stylized model for an exam schedule that minimizes the number of students with \emph{overlapping} exams, i.e., the number of students with two distinct exams assigned to the same time slot. 
\cite{Wel67} recognized the problem of finding an overlap-free exam schedule with as few time slots as possible as equivalent to a graph-coloring problem: given an undirected graph with a vertex for each class and edges exactly between any pair of classes that share students, finding an overlap-free schedule is equivalent to finding a minimum-vertex coloring.  Note that this formalization implies that even a simple version of exam scheduling that only aims to find a no-overlap schedule is in general NP-complete 
(see also \citeauthor{Gar79}~\citeyear{Gar79}).

Since then, work on exam scheduling  has spanned many directions.  
Foundational surveys by \cite{Car86}, \cite{Car95}, and \cite{Qu09} have covered the seminal work on graph-based and multi-criteria models, as well as on the use of techniques based on constraint programming; local search, including tabu search and simulated annealing; populations of solutions, including evolutionary algorithms, memetic algorithms, ant algorithms, and artificial immune algorithms; hyper-heuristics; and decomposition and clustering. 
Surveys with recent advances include \cite{Gas18},    \cite{Ald19}, \cite{Ces23}, and  \cite{Sie24}.  

Nevertheless, universities still face the same problem.  We conjecture that part of the reason for the lack of widespread adoption of exam-scheduling algorithms despite over six decades of active work is a divide that can be observed in the literature since at least the 1980's: a tension between (a) general variants of the exam scheduling problem, which may not be flexible enough to meet the practical constraints that vary from university to university; and (b) variants of the exam scheduling problem that are too specialized to a given institution, and thus difficult to apply elsewhere.  

General exam scheduling formulations and instances include the Uncapacitated Examination Timetabling Problem (UETP; see  \citeauthor{Lap84}~\citeyear{Lap84} and  \citeauthor{Car96}~\citeyear{Car96}) and the benchmark set from the 2nd International Timetabling Competition examination track at ITC 2007 (see \citeauthor{Mcc07}~\citeyear{Mcc07} and \citeauthor{Mcc12}~\citeyear{Mcc12}).  The ITC 2007 formulation is substantially more realistic than the UETP formulation, incorporating criteria such as minimizing the number of students with back-to-back exams and multiple exams in the same day, 
which are relevant to us. However, it is generally not an exact match for any particular institution.  For instance, it does not account for several of the criteria \BU considers.

When it comes to solving the problem for a specific institution,
there are a variety of formulations and approaches.
For example, \cite{Meh81} implemented a graph-based heuristic algorithm at Cedar Crest College to find exam schedules with few overlapping exams; \cite{Gar19} used integer programming to schedule assessments at the Universidad Polit{\'e}cnica de Madrid; and  \cite{Str17} used a variety of solvers with different modeling conventions in parallel 
to handle course, room, and exam scheduling at the Technion.  
See also \cite{Bac22}, \cite{Dim01}, \cite{Lap84}, \cite{Lot91}, \cite{Wan10}, and \cite{Ye24}, among many others.   
Of particular note, \cite{Ye24} involved work by a team of undergraduate student researchers who collaborated with the Cornell Registrar each semester and inspired our exam scheduling work involving undergraduate researchers at \BU.  
They provide a portfolio of schedules to their Registrar's Office by first assigning exams to exam blocks, then assigning exam blocks to time slots, and finally iteratively improving the resulting schedule. 
They use a layer-cake heuristic to schedule batches of courses one-at-a-time, with previously-scheduled courses warm-starting the next scheduling model, and allowing those previously scheduled courses to be rescheduled if needed.

While there has been broad work on exam scheduling, there has been limited work surveying what criteria matter most to stakeholders. 
\cite{Bur96b} surveyed the Registrars of 95 British universities specifically about exam scheduling tools. \citeauthor{Bur96b} asked them to rank the importance of 13 different constraints and criteria that a university might use in scheduling exams, 
and in turn the respondents suggested an additional 17.  The paper ends with a foreboding quote from one respondent about developing systems for exam scheduling: ``Having tried to do in the past what you are proposing now, I hope you have the time and resources to match your enthusiasm.  Do \underline{not} underestimate the complexity of the problem!  Good luck anyway.''   
 In the 2000's,  \cite{Cow02} surveyed an undisclosed number of students and invigilators (proctors), noting that students generally preferred to have gaps between exams and wanted a ``uniform distribution of exams over the examination period.'' Finally, in 2014,  \cite{Muk17} surveyed 50 students at the University of Nottingham.  Students conveyed a mixed opinion, with 74\% preferring an extended schedule and 12\% preferring a concentrated schedule of exams, 
 but they generally wanted to avoid having multiple exams on the same day.  

Finally, a small portion of the literature emphasizes designing tools that a Registrar's office can use directly. For instance, \cite{Gul20} implemented an integer programming model through a user-friendly spreadsheet for the Industrial Engineering Department of Yıldız Technical University;  \cite{Alh20} developed a graphical user interface for the German Jordanian University Registrar to use their CPLEX-based solver.

\section{Exam Scheduling at \BU}\label{sec:atBU}

\epigraph{\textit{If there's a way to design it to not have final exams, that would work best!}}{Comment from a respondent of our student survey on exams at \BU.}

In this section, we give context for final exam scheduling at \BU.  We begin by describing the general process, timeline, and  rules in Section \ref{sec:historical}.  Then we describe the results of a survey giving insight into how \BU students weigh tradeoffs in Section \ref{sec:survey}.

\subsection{\BU's Historical Exam Scheduling Process}\label{sec:historical}

At \BU, final exams have been historically scheduled by the registrar. \BU has nearly 4,000 undergraduate students and, in a typical semester, offers exams for between 500 and 600 courses.  Exams are scheduled in 22 three-hour time slots over 6 business days. 
The schedule has three daytime slots for each of those days.
Except for Friday and for the last day of exams, 
the schedule also has one nighttime slot each day. 
Before our collaboration,  the registrar would spend several months putting together an exam schedule to broadly support students as follows.

First, individual sections of a course are assigned to one of  approximately 70 \emph{course groups}, where all sections in a course group will be assigned the same exam slot. As a first step, a course group is created for each multi-section course that requests a common exam slot (e.g., a course group for all sections of Calculus 1).  Following that step, any courses not yet assigned to a course group are assigned to one corresponding to their regular meeting time (e.g., all classes that meet on Monday, Wednesday, and Friday (MWF) from 10:00 AM to 10:50 AM that are not already assigned a course group).  Since some courses have multiple types of meeting times (e.g., lecture times, a separate lab time, or an occasional evening time for midterms) there are niche policies followed by the registrar regarding which of such meeting times should be used to group the class. 

Second, the registrar assigns a time slot to each course group, with multiple course groups generally assigned to the same slot.  The registrar aims to minimize the number of students with overlapping exams and who have three or more exams fully contained within 24 hours, which are the cases in which a student can reschedule an exam. 
As noted before, it may be impossible to avoid requiring at least some students to have overlapping exams: 
students can sometimes concurrently enroll in two courses whose meeting times partially overlap and which are part of the same course group, and thus will always be assigned the same exam slot when following \BU policy. We ignore such cases in our results and report only \emph{unforced overlaps}, which do not account for those unavoidable situations.

In this process, the registrar also accounts for a variety of formal and informal criteria, such as avoiding students having back-to-back exams or a night exam followed by a morning exam on the next day; 
spacing out exams for students and especially for first-years; anticipating which exam slots are likely to have a conflict with athletics;  avoiding having a large number of students have to stay until the end of examination period; meeting faculty preferences for exam slots that support conference travel; etc.  (Note: throughout this paper, we use the term \emph{faculty} to describe instructional staff who teach classes and assign final exams.)

One challenge with the manual process comes from the difficulty in obtaining data for supporting it. 
In particular, the registrar had limited internal tools to access supporting data 
and looking up relevant data was tedious and time-consuming. 
For instance, counting overlaps and back-to-back exams required knowing the number of students concurrently enrolled in two courses.  Getting this information took multiple minutes for just a single pair of courses: the registrar would go to a website to get data reports, select a semester, wait for the semester data to load, select two class sections, wait for the report to be generated, see  a list of student IDs with that overlap, and then manually count the number of IDs.  This process only worked for individual course sections, and it did not allow the registrar to readily access data about overlaps between course groups.

Another challenge with the manual process comes from 
\BU's exam schedule having to be approved by three college deans before the semester begins.  This makes scheduling for the fall particularly tricky: full information about the incoming class of first-year students, such as their majors, AP scores, and course preferences,  is not known until late July, while the Fall semester begins mid-August.  Since the registrar's previous manual process took months, the registrar often drew on historical data to anticipate information about first-year students, and the registrar has described the exam scheduling process as akin to solving a Sudoku puzzle over several months.

\subsection{Student Preference Data}\label{sec:survey}

Before modeling the problem, we surveyed students at \BU to better understand their preferences regarding final exam schedules.  Students were asked to rate nine different possible exam inconveniences.  The survey was open to all students, of whom 226 responded, and their responses and the exam inconveniences are shown in Table \ref{fig:studentResp}.  Students generally were in alignment that congested exam schedules were the most inconvenient, but far preferred having exams on the first day to the last day of the exam period. This may surprise faculty since, at \BU, students only have one reading day and no weekend between when classes end and when exams begin. Students across our different programs were generally in agreement, except on having an exam on the last day: 
only 25.6\% of the engineering students rated that criteria as ``Extremely inconvenient,'' compared to 38.2\% of the students in arts and sciences and 45.05\% of the management students.

\begin{table}[t]
    \centering
    \includegraphics[width=0.8\linewidth, trim={2.4cm 17.4cm 2.4cm 3cm}, clip]{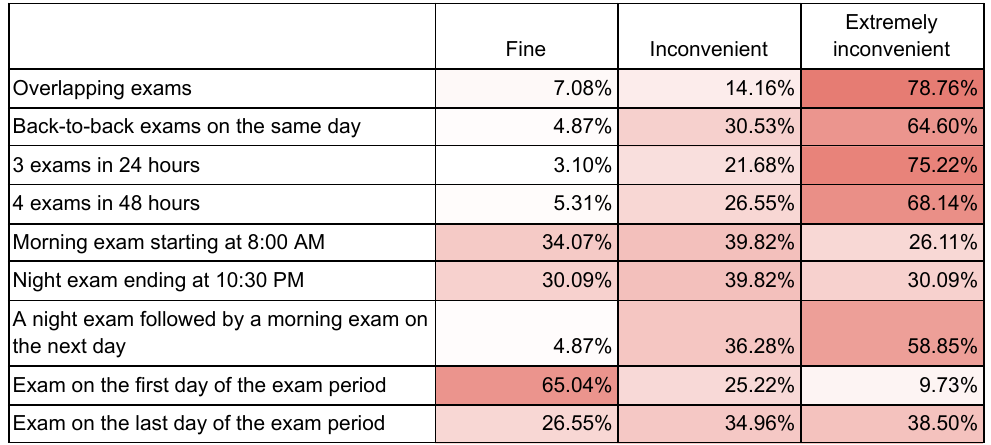}
    \caption{Student responses to final exam schedule survey. Students were asked to rate each of nine exam inconveniences.}
    \label{fig:studentResp}
\end{table}

Students were also given an open-ended question: \emph{Are there any other factors we should consider when designing the final exam schedule, either that make it better or worse?}  Unsurprisingly, student responses revealed different preferences, ranging from  ``No 8am exams please'' to ``no night exams.'' There were also impractical suggestions, e.g., ``If there's a way to design it to not have final exams, that would work best!''  Several students suggested trying to incorporate the subjective difficulty of the course, e.g., ``How hard the course is. Preferably the `harder' courses should have more time to study!''

Based on the results of the survey, we have focused on avoiding the five inconveniences that concerned students the most: 
overlapping exams, 3 exams in 24 hours, 4 exams in 48 hours, back-to-back exams on the same day, and night-to-morning exams. 
We use a flexible weighting scheme for these factors, and then generate a portfolio of exam schedules with different weightings; one of the weighting schemes is based on the relative preferences of students from our survey.

\section{A Collaborative Approach to Exam Scheduling}\label{sec:Approach}

In this section, we discuss our approach to exam scheduling at \BU, which involved close collaboration with the Registrar's Office. 
We began preliminary work in Spring 2023 with a series of meetings to understand their current process as described in Section~\ref{sec:historical}, 
as well as by obtaining IRB approval and conducting the survey with students as described in Section~\ref{sec:survey}. 

We present subsequent steps of our collaboration in the subsections below.
First, in Summer 2023, we began working in earnest and developed and released visualization tools that helped the registrar with their historical manual scheduling process (Section~\ref{sec:vis}).
Second, we started developing integer programming models and heuristics for scheduling exams in Summer 2023 and continued through Summer 2024.  During this time, we worked closely with the registrar to meet additional requirements and improve the time to obtain solutions.  
We initially relied on the academic license of a commercial solver, 
and then adapted our model and algorithms to work with an open-source solver (Section~\ref{sec:IP}).
Finally, through Fall 2024, we designed and refined an open-source user interface that allowed the registrar to directly build and adapt a portfolio of generated schedules based to their needs (Section~\ref{sec:tool}).  The registrar first used our visualization tools to schedule the Fall 2023 exams, and began using our optimization tools and interface to schedule the Spring 2024 exams.

\subsection{Visualization Tools}\label{sec:vis}

To understand and support the registrar's needs, 
we began by building tools to expedite their historical process.  
First, 
we gathered the anonymized enrollment information that would be necessary for evaluating exam schedules, such as student enrollment information and course meeting times.   
Then we developed Tableau dashboards to provide summary information relevant to exam scheduling -- such as the number of students mutually enrolled in two course groups --  for both the current and recent semesters.  
The registrar used this dashboard  to help manually schedule the Fall 2023 final exams, 
before we  provided automated schedules.

\begin{figure}
    \centering
    \includegraphics[width=0.98\linewidth, trim=0 0.45cm 0 0, clip]{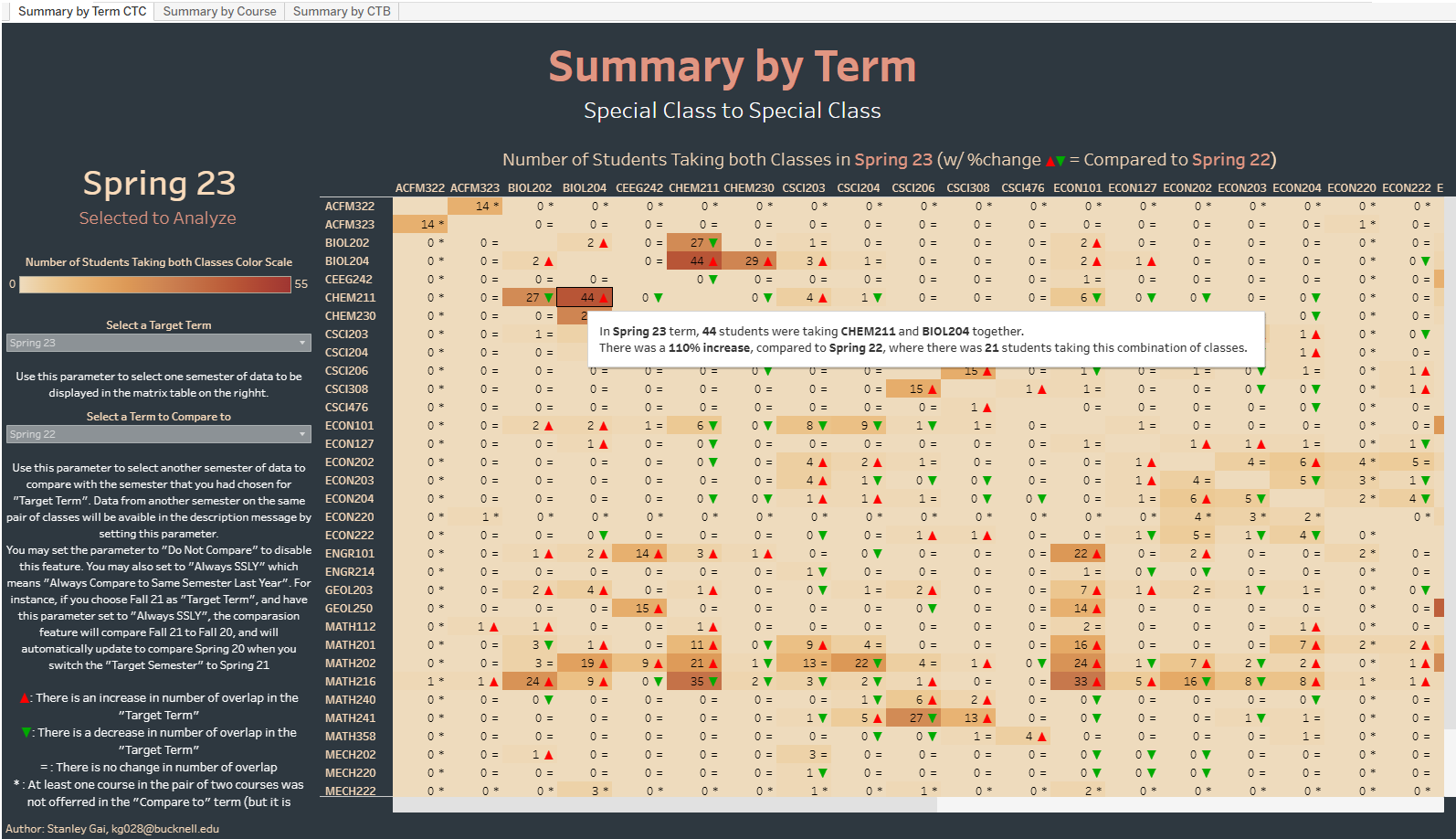}
    \caption{Tableau interface for visualizing the overlap information between pairs of course groups, showing both the overlap for the current semester and how it changed from a previous academic year.  Hovering over a specific pair of course groups displays more detailed information.}
    \label{fig:tab1}
\end{figure}

\begin{figure}
    \centering
    \includegraphics[width=0.98\linewidth]{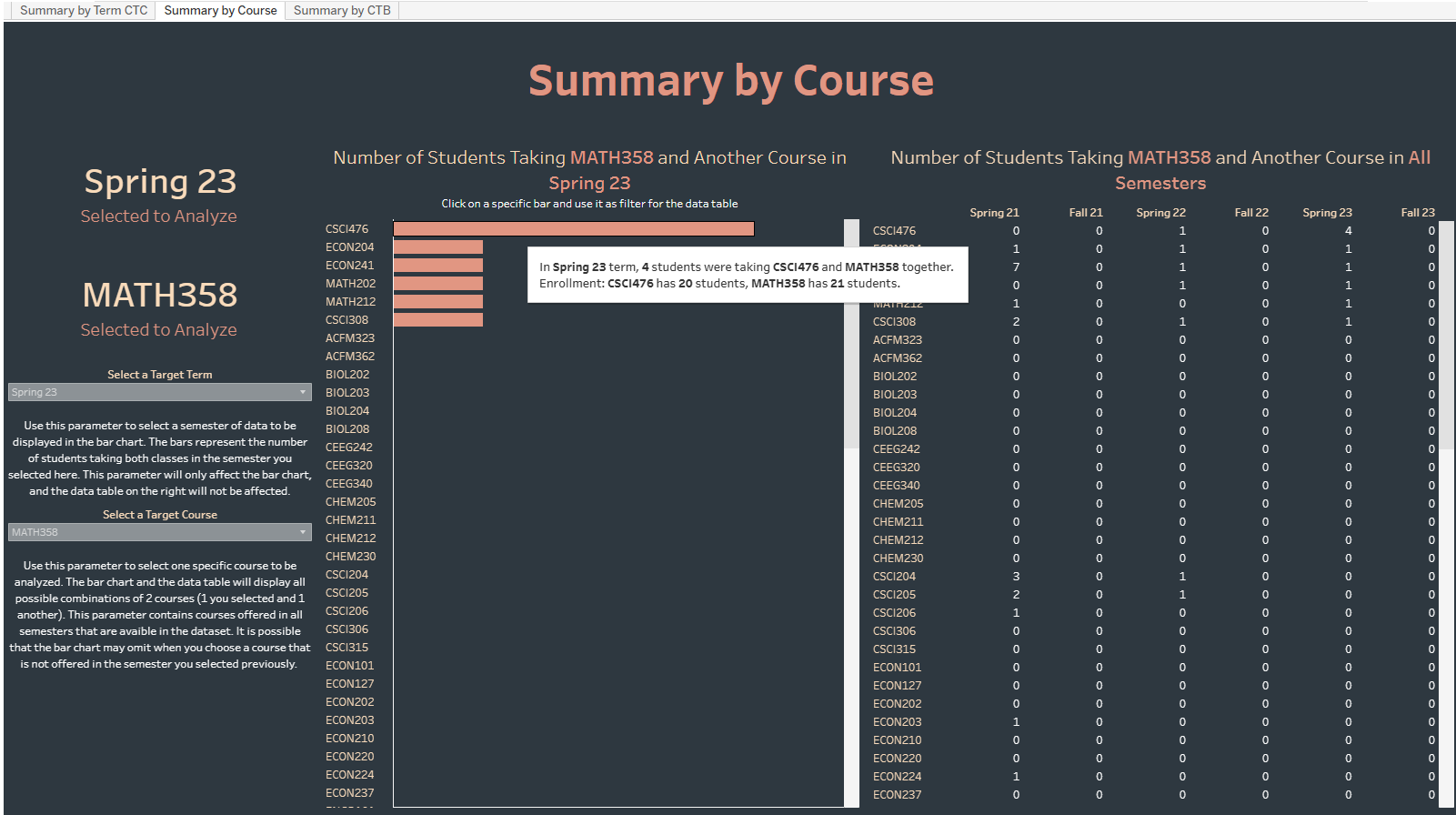}
    \caption{Tableau interface for visualizing overlap information about a specific course group, highlighting top conflicts from the current semester and providing detailed information about historical overlaps.}
    \label{fig:tab2}
\end{figure}

Figure \ref{fig:tab1} shows the first pane of the dashboard.  The user selects a current term (Spring 2023 in the figure), and a grid shows the number of students concurrently enrolled in any pair of course groups (on this pane, specifically course groups based on classes requesting common exam times).  For instance, the highlighted cell shows that 44 students were concurrently enrolled in any section of BIOL204 and in any section of CHEM211.  The user can also select a previous term to view historical trends (Spring 2022 in the figure).  Red and green arrows indicate if there is a change in the number of concurrent enrollments from the previous academic year, and hovering over a cell shows more detailed information.  In Figure~\ref{fig:tab1}, the hover text shows that only 21 students were concurrently enrolled in BIOL204 and CHEM211 in the previous term.  We included this historical data since the registrar must schedule exams before add-and-drop begins and while enrollments are still occurring.  The dashboard also allows the registrar to visualize data about a specific course (say, when deciding where to manually place a course)  as in Figure \ref{fig:tab2}, as well as data for course groups based on meeting time in another pane.

\subsection{Mathematical Formulation and Heuristic Solution}\label{sec:IP}

In parallel to developing the dashboards, we began formulating an integer programming model for scheduling exams.  We initially formulated a model which enforced several constraints:
\begin{itemize}
    \item Every course group must be scheduled for an exam in exactly one time slot.
    \item The number of distinct students with an exam in each time slot is capped based on classroom capacity.
    \item Any manually-entered requirements for certain course groups, such as assigning them to or blocking them off from specific time slots, must be satisfied.
\end{itemize}
Subject to meeting those constraints, the model aimed at minimizing a weighted combination of the student inconveniences listed in Section~\ref{sec:survey} (overlapping exams, 3 exams in 24 hours, 4 exams in 48 hours, back-to-back exams on the same day, and night-to-morning exams) as well as two similar forms of inconvenience for faculty (overlapping exams and back-to-back exams on the same day). 
Changing the relative weights of those inconveniences allowed us to produce a portfolio of schedules for the registrar to use as starting points. 
Through the interface described in the next section, 
each of those schedules could be directly modified by the registrar.   We refer the reader to Appendix~\ref{sec:ipf} for a description of the integer programming formulation. 

While the exact number of constraints and variables changes every semester, 
the resulting model had approximately 300,000 to 400,000 binary decision variables and 300,000 to 400,000 constraints. 
This model would generally not solve to optimality, 
and the solutions obtained did not always compare favorably to the manual solution. 
Hence, we implemented a heuristic that breaks the problem in two phases. 
First, in phase 1, we only schedule exams for large-enrollment course groups and the small handful of courses that must be assigned to certain days or times, and we only consider some of the inconveniences. 
Then, in phase 2, we solve the entire problem with additional constraints fixing the exam slots of the course groups from phase 1 according to the solution obtained for it.

\begin{figure}
    \centering
  \includegraphics[width=.9\linewidth]{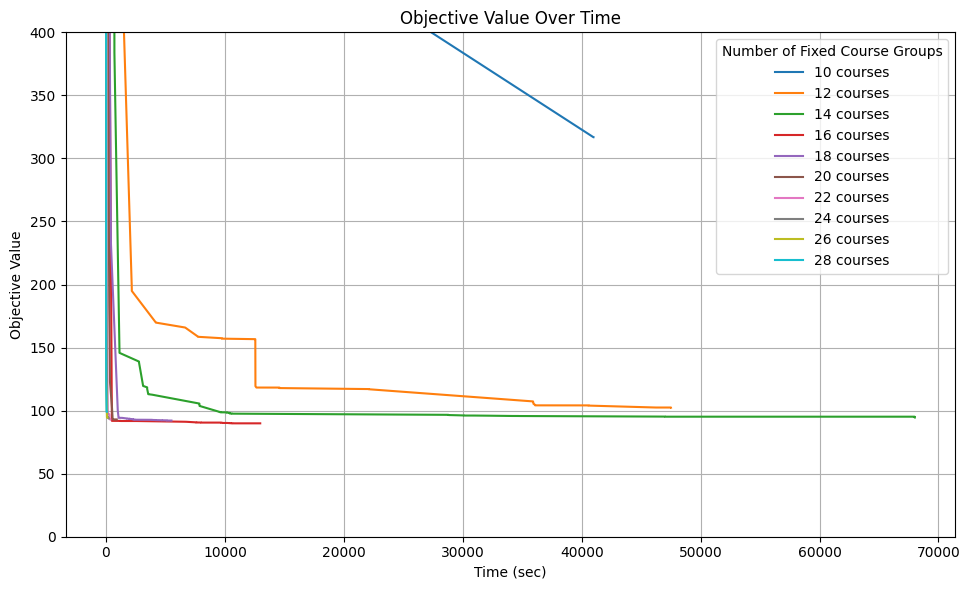}
    \caption{Final phase 2 objective value versus phase 2 running time based on the number of courses fixed in phase 1 for a sample semester.}
    \label{fig:numcourses}
\end{figure}

\begin{figure}
    \centering
  \includegraphics[width=.9\linewidth]{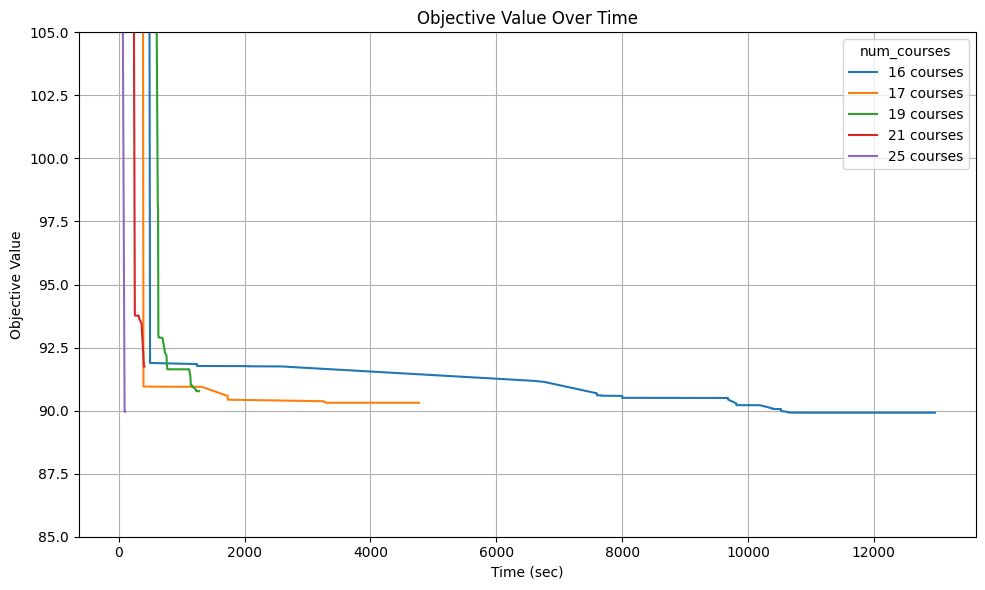}
    \caption{Zoomed in view of phase 2 objective value versus phase 2 running time based on the number of courses fixed in phase 1 for a sample semester, only including the 5 best runs.}
    \label{fig:numcourses2}
\end{figure}

Figures \ref{fig:numcourses} and \ref{fig:numcourses2} show sample runs with different numbers of fixed course groups.  Phase 2's running time was capped at 24 hours, but all trials with at least 16 fixed course groups ran to completion. Note that there is substantial variance in both how long it took for phase 2 to run and the quality of the final solution. For the semester shown in Figures \ref{fig:numcourses} and \ref{fig:numcourses2}, it took under a minute for phase 2 to terminate with 28 fixed course groups, while it took over 18 hours for phase 2 to terminate with 14 fixed course groups. Generally speaking, fewer fixed course groups required more time in phase 2.  Fixing 10 course groups led to a solution that was more than three times as bad as any other run.  In general, over many semesters of testing, we found that fixing anywhere from 17 to 21 course groups tended to find the best solutions in reasonable amounts of time, and that with those choices, optimal phase 2 solutions were generally found within 2 hours. 

Even though phase 2 is roughly the same size as our original formulation, the constraints fixing a small number of course groups break the inherent symmetry in exam scheduling and allow the phase 2 model to generally be solved to optimality. As an indicative example using the open-source solver SCIP and fixing 20 course groups, phase 2 initially had a model 352,056 variables and 306,669 constraints; after presolve, the model only had 76,722 variables and 117,970 constraints -- of which 58,262 were added during presolve.  With 24 course groups fixed, phase 2 only had 8,387 variables and 16,458 constraints after presolve. 

However, even if phase 1 and phase 2 both solve to optimality, we are not guaranteed a globally optimal solution. Thus, in practice we generate a portfolio of schedules and set reasonable time limits on each phase; for phase 1, we keep extending the time limit as we find new incumbents and for phase 2 we set a maximum runtime of 4 hours.  Our tool produces 20 possible exam schedules: We consider four different sets of inconvenience weights and, and for each set of inconvenience weights, generate five schedules where the number of course groups fixed in phase 1 ranges from 17 to 21.  The tool ultimately reports four possible schedules to the user: the best schedule found for each set of inconvenience weights.  Our server allows four schedules to be generated in parallel, and the registrar can generally generate a complete portfolio overnight.  We refer the reader to Appendix~\ref{sec:heu} for additional information about the two-phase heuristic.
 
Table \ref{tab:res} shows our heuristic being evaluated on historical data, using the semesters immediately before we began working with the Registrar's Office.  These results show substantial improvement: Far fewer students have overlapping exams and three exams within 24 hours. These are the most consequential inconveniences at \BU, since they require moving one exam from its scheduled slot for that student.  Similarly,  there is a substantial reduction in the number of students with back-to-back exams, the most frequent inconvenience.

\begin{table}[ht]
\centering
\tableFont
\caption{Evaluating our heuristic on historical data:  Real shows the number of inconveniences in the manually-produced schedule for that semester, Opt shows the number of inconveniences using our heuristic, and Imp shows the percent Opt/Real rounded to the nearest percent.}\label{tab:res}
\begin{tabular}{@{}l ccc ccc ccc@{}}
\toprule
\textbf{} & \multicolumn{3}{c}{\textbf{Spring 2022}} &  \multicolumn{3}{c}{\textbf{Fall 2022}} & \multicolumn{3}{c}{\textbf{Spring 2023}} \\
 & Real & Opt & Imp & Real & Opt & Imp & Real & Opt & Imp \\
\midrule
Students with an Unforced Overlap & 37 & 1 & 3\% & 44 & 0 & 0\% & 37 & 0 & 0\% \\
Students with 3 Exams in 24 Hours & 155 & 74 & 48\% & 118 & 45 & 38\% & 229 & 74 & 32\% \\
Students with 4 Exams in 48 Hours & 67 & 31 & 46\% & 55 & 22  & 40\% & 75 & 36 & 48\% \\
Students with Back-to-Back Exams & 701 & 381 & 54\% & 668 & 354 & 53\% & 1085 & 423 & 39\% \\
Students with Night-to-Morning Exams & 185 & 86 & 46\% & 136& 63 & 46\%& 198 & 75 & 38\% \\
Students with at Least One Inconvenience & 849 & 450 & 47\% & 798 & 400 & 50\% & 1166 & 488 & 42\% \\
Faculty with an Unforced Overlap & 12 & 10 & 83\% & 9 & 9& 100\%& 12 & 12 & 100\% \\
Faculty with Back-to-Back Exams & 17 & 4 & 24\% & 5 & 7 & 140\% & 11 & 9 & 82\% \\
\bottomrule
\end{tabular}
\end{table}

\vspace{1em}

Finally, we note that our implementation used for these experiments was fully open-source, using the solver SCIP.  We compared implementations in both Gurobi and SCIP. Table \ref{tab:scip_vs_gurobi} shows how the runtime and final solution found compares across four semesters.  In three of the four semesters, Gurobi obtained a slightly better final solution, whereas SCIP obtained a better solution once. However, the final objective values were always within 4\% of each other.  In each semester, however, Gurobi was substantially faster: phase 2 generally took minutes to run with Gurobi versus 1 to 2 hours with SCIP.

\begin{table}[ht]
\centering
\tableFont
\caption{Comparison of SCIP to Gurobi for Objective Value and Running Time. Percentages represent SCIP as a proportion of Gurobi.}\label{tab:scip_vs_gurobi}
\begin{tabular}{@{}l cccc@{}}
\toprule
\textbf{} & \textbf{Spring 2022} & \textbf{Fall 2022} & \textbf{Spring 2023} & \textbf{Fall 2023} \\
\midrule
Objective Value: SCIP/Gurobi & 103.58\% & 101.51\% & 99.91\% & 100.43\% \\
Running Time: SCIP/Gurobi & 1320.78\% & 1458.72\% & 3271.22\% & 1428.42\% \\
\bottomrule
\end{tabular}
\end{table}

\subsection{An Open-Source Web-Based User Interface}\label{sec:tool}

As the final step of our process, we designed and built an open-source tool that allows the registrar to build a portfolio of optimized final exam schedules themselves.  This tool runs the two-phase heuristic that we developed in an easy-to-use web-browser-based interface and  does not require any mathematical or technical background to use.  The Registrar's Office used this tool to schedule the Fall 2024 and Spring 2025 exams, while we fine-tuned it based on their feedback.  

\begin{figure}
    \centering
  \includegraphics[width=.7\linewidth]{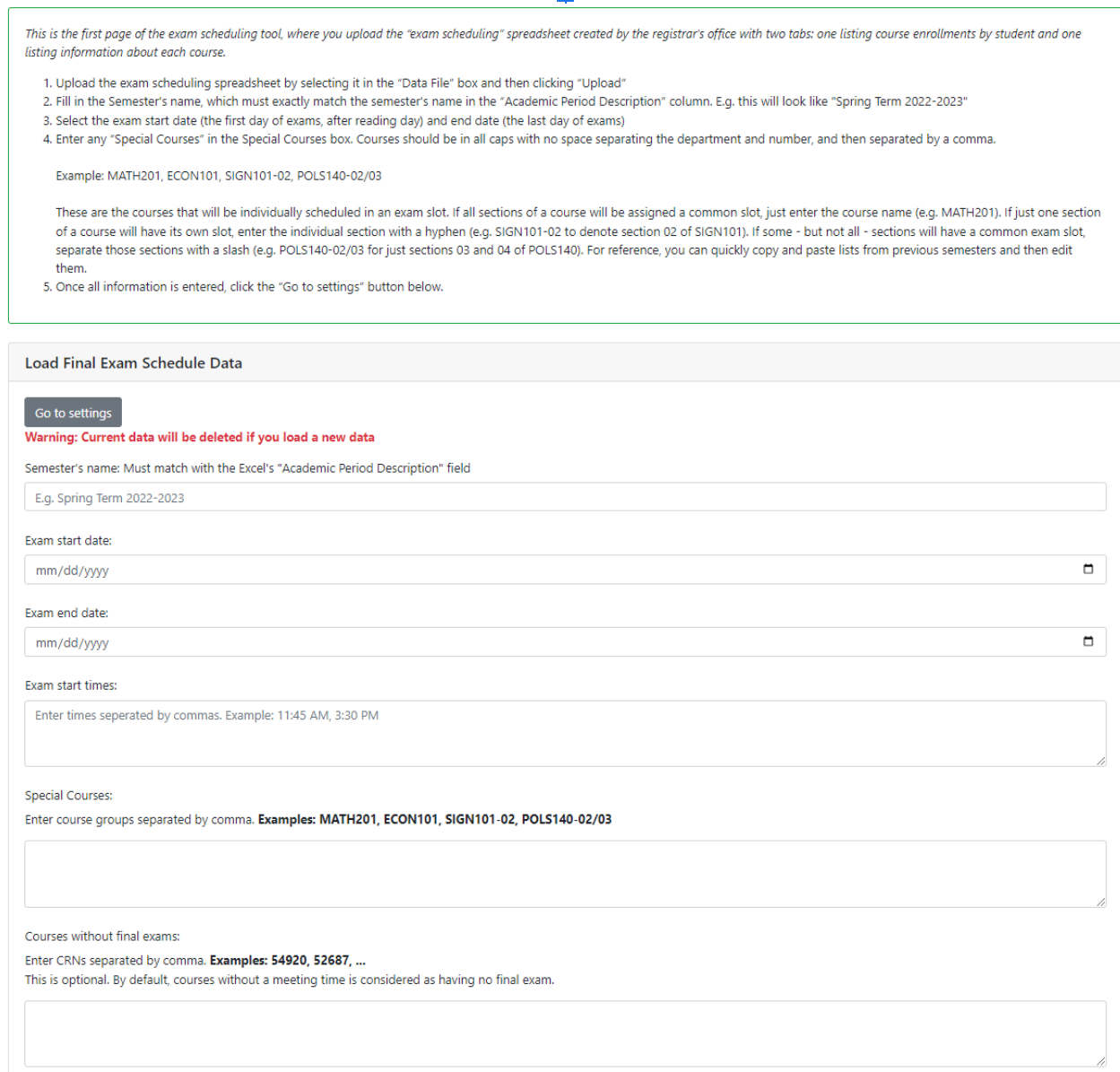}
    \caption{Interface for entering initial information for exam scheduling.}
    \label{fig:sem}
\end{figure}

\begin{figure}
    \centering
  \includegraphics[width=.7\linewidth, trim=0 0 0 5.2cm, clip]{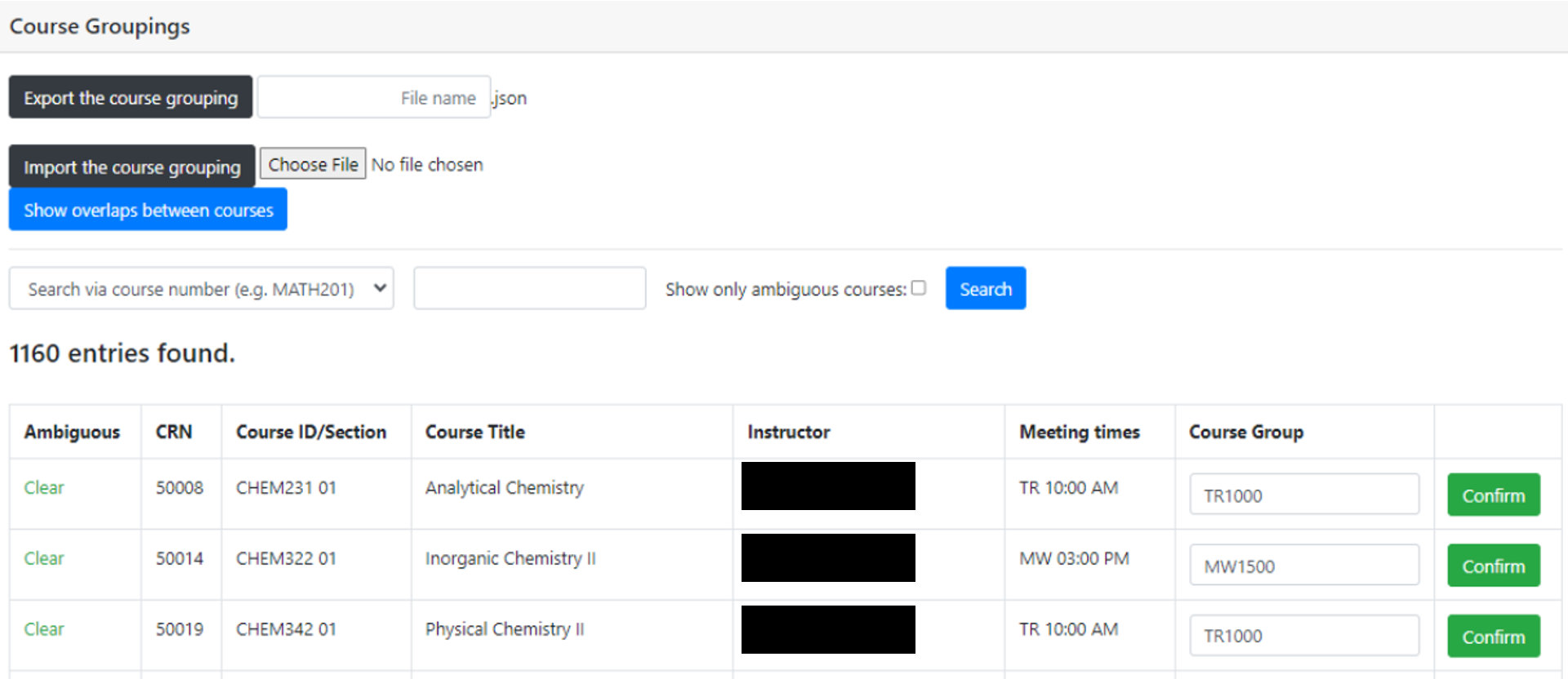}
    \caption{Interface for confirming and editing course groups.}
    \label{fig:group}
\end{figure}

\begin{figure}
    \centering
  \includegraphics[width=.95\linewidth, trim=0 9.9cm 0 0, clip]{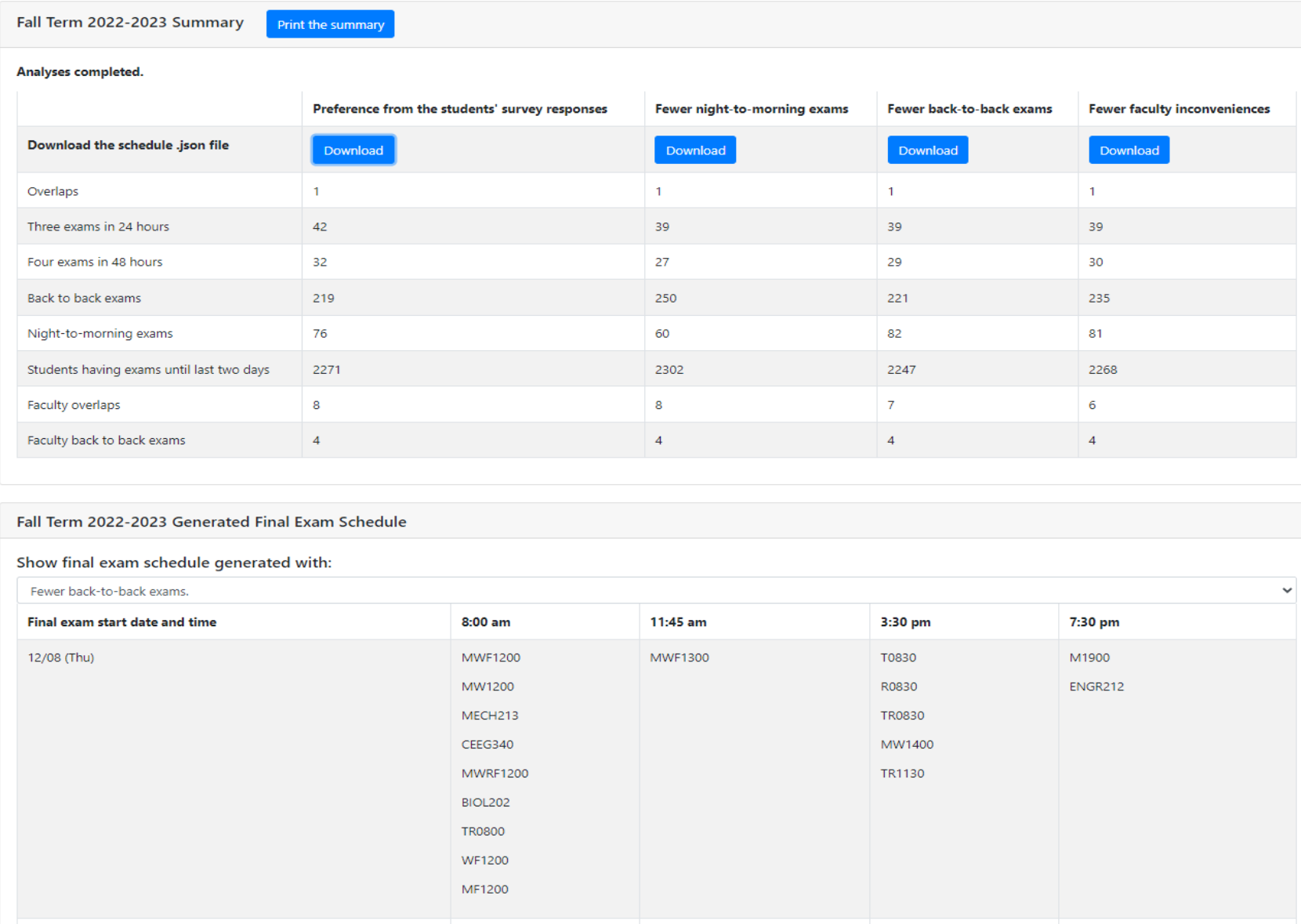}
    \caption{Interface to access a portfolio of optimized schedules.}
    \label{fig:port}
\end{figure}

\begin{figure}
    \centering
  \includegraphics[width=.8\linewidth]{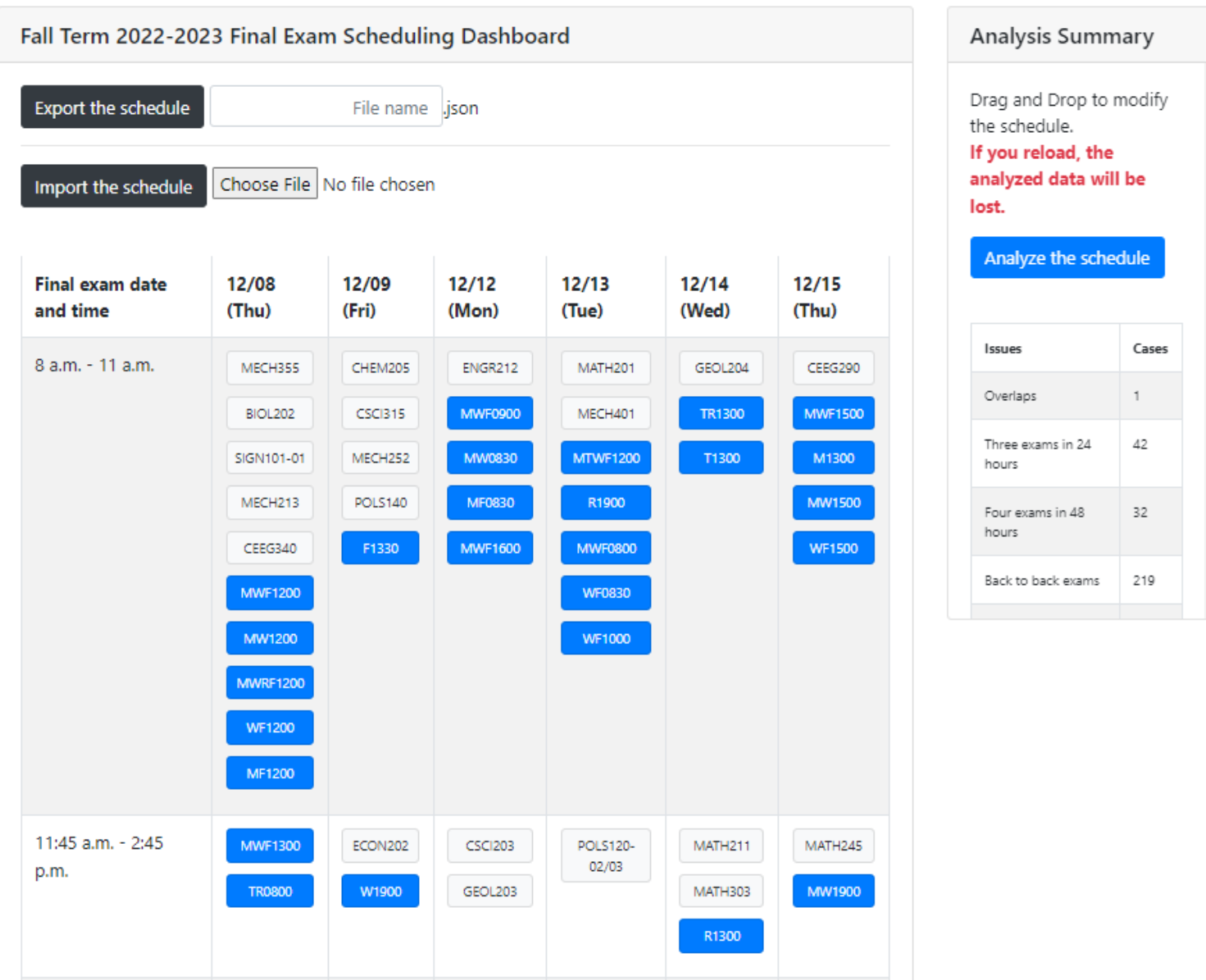}
    \caption{Interface for viewing and modifying a schedule with a ``drag-and-drop'' feature.}
    \label{fig:view}
\end{figure}

Figures \ref{fig:sem} through \ref{fig:view} show features of this tool.  The registrar opens the tool through a web-browser, and the tool is automatically linked to enrollment data.  When creating an exam schedule, the registrar first provides core information about that semester, such as listing any multi-section courses that have requested a coordinated exam; see Figure \ref{fig:sem}.  Courses are then automatically grouped based on that list and their meeting times, but these course groupings can be checked and modified by the registrar; see Figure \ref{fig:group}.  Courses with multiple meeting times that do not unambiguously correspond to a course group are automatically flagged as ambiguous for the registrar to check.  Next, the registrar generates a portfolio of schedules.  They can enter additional constraints, such as restricting certain exams to certain times, and then see a summary of the portfolio of optimized schedules as in Figure \ref{fig:port}.  Finally, the registrar can make additional changes to a schedule through a ``drag and drop'' interface, directly seeing how those changes affect student and faculty inconveniences; see Figure \ref{fig:view}.

We note that this tool is fully open-source and tailorable to a variety of university settings: it captures common notions of student inconveniences whose relative weights can be easily modified to meet a school's priorities; it is flexible for schools with different numbers of exam slots; and it only requires limited data about student enrollment and faculty course assignments.  Our tool is available on \GIT, including a detailed guide for how to use the tool written for a non-technical audience and a README for customizing the code for a more technical audience.

\section{Results and Impact}\label{sec:Results}

The Registrar's Office manually created final exam schedules using the historical process outlined in Section \ref{sec:historical} through the Spring 2023 semester.  In Fall 2023, they used our Tableau dashboard to facilitate their manual process.  Beginning in Spring 2024, schedules were generated using our tool and then finalized by the Registrar's Office.

\begin{figure}
    \centering
  \includegraphics[width=.99\linewidth]{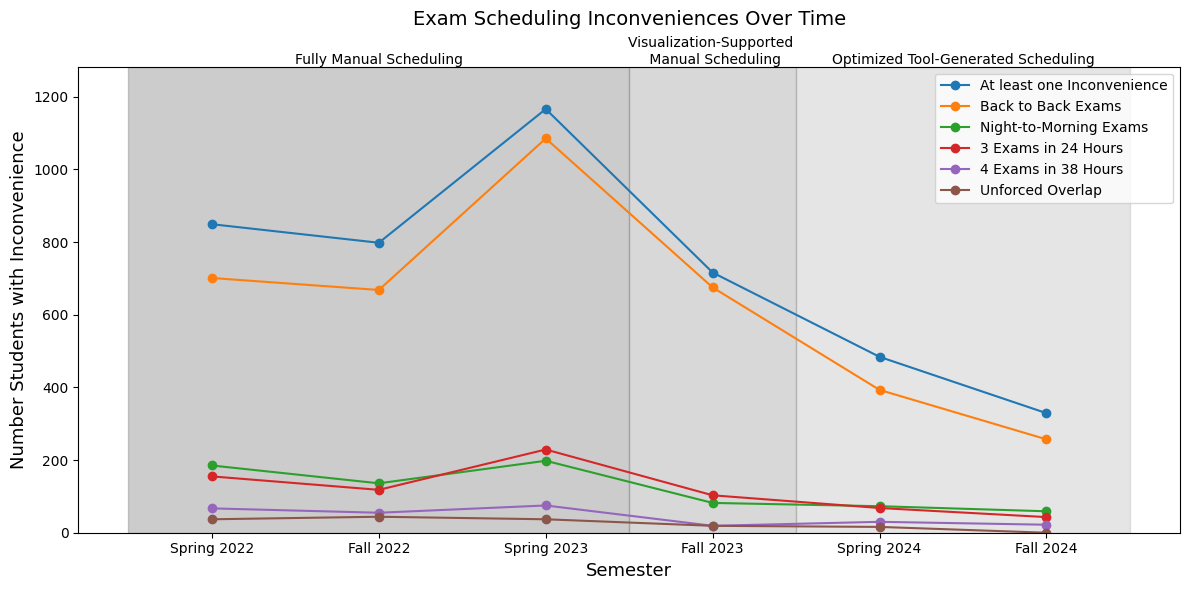}
    \caption{The numbers of students with final exam inconveniences before and since using our tools.}
    \label{fig:eval}
\end{figure}

Combining the \BU registrar's insight into exam scheduling with our tool has led to notably better exam schedules.  Table \ref{tab:eval} summarizes attributes of exams scheduled before this project, with the support of the Tableau dashboard in Fall 2023, and with our optimization approach in Spring 2024 and Fall 2024;  Figure \ref{fig:eval} shows the same data visually.  Note that we have seen major reductions in the number of students with unforced overlaps and three exams in 24 hours, as well as in the number of students with back-to-back exams and with at least one inconvenience.

\begin{table}[ht]
\centering
\tableFont
\caption{The numbers of faculty and students with final exam inconveniences before and since using our tools. (For conciseness, S stands for Spring and F for Fall.)}\label{tab:eval}
\begin{tabular}{@{}l ccc c cc@{}}
\toprule
\textbf{} & \multicolumn{3}{c}{\textbf{\makecell{Fully Manual\\Scheduling}}} & \textbf{\makecell{Visualization-\\Supported\\Manual\\Scheduling}} & \multicolumn{2}{c}{\textbf{\makecell{Optimized\\Tool-Generated\\Scheduling}}}\\
 & S 2022 & F 2022 & S 2023 & F 2023 & S 2024 & F 2024 \\
\midrule
Students with an Unforced Overlap & 37 & 44 & 37 & 19 & 16 &0 \\
Students with 3 Exams in 24 Hours & 155 & 118  & 229 & 103 & 68& 43  \\
Students with 4 Exams in 48 Hours & 67 & 55 & 75 & 19  & 30 & 22 \\
Students with Back to Back Exams & 701 & 668 & 1085 & 675 & 393 & 257 \\
Students with Night-to-Morning Exams & 185 & 136 & 198 & 82 & 73 & 59 \\
Students with at Least One Inconvenience & 849 & 798 & 1166 & 716 & 484 & 329 \\
Faculty with an Unforced Overlap & 12 & 9 & 12 & 8  &9 & 1 \\
Faculty with Back to Back Exams & 17 & 5 & 11 & 9 & 15 & 5 \\
\bottomrule
\end{tabular}
\end{table}

Finally, we note that this tool also allows us to explore ``what-if'' situations.  For instance, Table \ref{tab:exam_day_comparison} explores what effect we would have by extending or shortening the final exam period by a day.  While the number of unforced overlaps is unchanged, removing an exam day substantially increases the number of students with three exams in 24 hours. Those, in turn, would need to have an exam rescheduled.  In contrast, adding a single exam day would roughly halve that number, as would it halve the number of students with at least one inconvenience.  

\begin{table}[ht]
\centering
\tableFont
\caption{What-if analysis showing the effect of increasing or decreasing the length of the final exam period (computed based on Spring 2025 data before the add-and-drop period).}\label{tab:exam_day_comparison}
\begin{tabular}{@{}l ccc@{}}
\toprule
\textbf{} & \textbf{\makecell{6 Exam Days\\(Current Schedule)}} & \textbf{5 Exam Days} & \textbf{7 Exam Days} \\
\midrule
Students with an Unforced Overlap & 10 & 10 & 10 \\
Students with 3 Exams in 24 Hours & 81 & 189 & 40 \\
Students with 4 Exams in 48 Hours & 76 & 89 & 44 \\
Students with Back to Back Exams & 419 & 763 & 220 \\
Students with Night-to-Morning Exams & 196 & 206 & 53 \\
Students with at Least One Inconvenience & 603 & 932 & 306 \\
Faculty with an Unforced Overlap & 8 & 10 & 9 \\
Faculty with Back to Back Exams & 7 & 11 & 2 \\
\bottomrule
\end{tabular}
\end{table}

\section{Conclusion}\label{sec:Conclusion}
This paper presented a study case, algorithms, and open-source software for final exam scheduling.  This work has been implemented at \BU, leading to an exam schedule that is much more conducive to student success.  Between  Spring 2022 and Spring 2023, which correspond to the period immediately before this project was carried out, an average of 39.33 students per semester had overlapping exams, and an average of 167.33 students had three exams in 24 hours. Both of those situations require an exam to be moved. In  Spring 2024 and Fall 2024, which correspond to the first semesters in which our tools were completely is use, those numbers dropped to 8 and 55.5, respectively.  Similarly, the average number of students with at least one inconvenience has more than halved, from 937.66 to 406.5.   As a result, this work has received national attention and been featured in \CHE in a series on data-based decisions to foster student success at universities.

Moreover, our tools were built in direct collaboration with the Registrar's Office.  By closely working with the registrar throughout the entire process, we were able to tailor tools that serve their actual needs; through collaboration, we were able to combine the registrar's domain expertise with analytics tools.  At the same time, this tool has saved substantial time for the Registrar's Office: a process that used to take months of manual work can now be done in a matter of days. 

Finally, while our specific implementation is based on \BU, we have designed the tools to be adaptable by other universities: they can be readily customized, rely fully on open-source optimization software, and account for a variety of common hard and soft constraints on exam schedules.

\beginFEappendix

\section{Integer Programming Formulation}\label{sec:ipf}

\noindent \textbf{Indices and Sets}
\begin{description}[style=multiline,leftmargin=2.5cm,font=\normalfont,noitemsep,topsep=0pt]
\item[$s \in \mathcal{S}$] The set of students who need to take final exams.
\item[$g \in \mathcal{G}$] The set of all course groups (e.g., Calculus I, or courses that meet MWF at 10:00 AM).
\item[$f \in \mathcal{F}$] The set of faculty members.
\item[$t \in \mathcal{T}$] The set of all time slots in which final exams could be scheduled.
\item[$\mathcal{T}_\text{B2B}$] The set of pairs of back-to-back time slots $\{t_0,t_1\}\subset\mathcal{T},$ meaning that $t_0$ immediately precedes $t_1$ on the same day.
\item[$\mathcal{T}_\text{PMtoAM}$] The set of pairs of time slots $\{t_0,t_1\}\subset\mathcal{T}$ where $t_0$ is a night time slot and $t_1$ is a morning time slot on the following day.
\item[$\mathcal{T}_{3in24}$] The set of subsets $T\subset\mathcal{T}$ having $|T|=3$ that consist of three time slots that fall within a 24-hour time period.
\item[$\mathcal{T}_{4in48}$] The set of subsets $T\subset\mathcal{T}$ having $|T|=4$ that consist of four time slots that fall within a 48-hour time period.\\
\end{description}

\noindent \textbf{Data}
\begin{description}[style=multiline,leftmargin=2.5cm,font=\normalfont,noitemsep,topsep=0pt]
\item[$a_t$] An indicator parameter for $t\in \mathcal{T}$. If $a_t = 1$, then time slot $t$ is available for final exams.
\item[$b_{sg}$] An indicator parameter for $s\in \mathcal{S}, g\in \mathcal{G}$. If $b_{sg} = 1$, then student $s$ is enrolled in a course that belongs to course group $g.$
\item[$d_{fg}$] An indicator parameter for $f\in \mathcal{F}, g\in \mathcal{G}$. If $d_{fg} = 1$, then faculty member $f$ is teaching a course that belongs to course group $g.$
\item[$r_{gt}$] An indicator parameter for $g\in \mathcal{G}, t\in \mathcal{T}$. If $r_{gt} = 1$, then course group $g$ is required to be assigned to time slot $t.$ 
\item[$q_{gt}$] An indicator parameter for $g\in \mathcal{G}, t\in \mathcal{T}$. If $q_{gt} = 1$, then course group $g$ is forbidden to be assigned to time slot $t.$ 
\item[$N_g$] The total number of students enrolled in course group $g \in \mathcal{G}.$
\item[$M_1$] The maximum number of students that can be assigned to take a final exam during any single time slot (due to space limitations, for example).
\item[$M_2$] The maximum number of courses that a student could take during a semester.
\item[$M_3$] The maximum number of courses that a faculty member could teach during a semester.
\item[$\rho_\text{overlap}$] A user-defined penalty parameter that is incurred whenever a student is enrolled in two course groups that are assigned the same time slot.
\item[$\rho_\text{B2B}$] A user-defined penalty parameter that is incurred whenever a student is scheduled to take two exams in back-to-back time slots on the same day.
\item[$\rho_\text{PMtoAM}$] A user-defined penalty parameter that is incurred whenever a student is scheduled to take a night exam followed by a morning exam the next day.
\item[$\rho_\text{3in24}$] A user-defined penalty parameter that is incurred whenever a student is scheduled to take three exams within a 24-hour time period.
\item[$\rho_\text{4in48}$] A user-defined penalty parameter that is incurred whenever a student is scheduled to take four exams within a 48-hour time period.
\item[$\rho_\text{facoverlap}$] A user-defined penalty parameter that is incurred whenever a faculty member is teaching courses in two course groups that are assigned the same time slot.
\item[$\rho_\text{facB2B}$] A user-defined penalty parameter that is incurred whenever a faculty member is scheduled to give two exams in back-to-back time slots on the same day.\\
\end{description}

\noindent \textbf{Binary Decision Variables}
\begin{description}[style=multiline,leftmargin=2.5cm,font=\normalfont,noitemsep,topsep=0pt]
\item[$x_{gt}$] 1 if course group $g$ is assigned to time slot $t$; 0 otherwise.
\item[$v_{st}$] 1 if student $s$ is scheduled to take an exam during time slot $t$; 0 otherwise.
\item[$w_{ft}$] 1 if faculty member $f$ is scheduled to give an exam during time slot $t$; 0 otherwise.
\item[$z_{st}^\text{overlap}$] 1 if student $s$ is enrolled in two course groups that are assigned to the same time slot $t$; 0 otherwise.
\item[$z_{st}^\text{B2B}$] 1 if student $s$ is scheduled to take two exams in back-to-back time slots, beginning with time slot $t$; 0 otherwise.
\item[$z_{st}^\text{PMtoAM}$] 1 if student $s$ is scheduled to take a night exam in time slot $t,$ followed by a morning exam the next day; 0 otherwise.
\item[$z_{s}^\text{3in24}$] 1 if student $s$ is scheduled to take three exams within a 24-hour time period; 0 otherwise.
\item[$z_{s}^\text{4in48}$] 1 if student $s$ is scheduled to take four exams within a 48-hour time period; 0 otherwise.
\item[$z_{f}^\text{facoverlap}$] 1 if faculty member $f$ is enrolled in two course groups that are assigned to the same time slot; 0 otherwise.
\item[$z_{f}^\text{facB2B}$] 1 if faculty member $f$ is scheduled to give two exams in back-to-back time slots; 0 otherwise.\\
\end{description}

\noindent \textbf{Constraints}
\begin{itemize}
    \item Each course group must be assigned to exactly one time slot.
    \begin{equation} 
        \sum_{t \in \mathcal{T}} x_{gt} =1 \quad \mbox{ for } g \in \mathcal{G} \label{const1}
    \end{equation}
    
    \item These constraints relate the $x_{gt}$ and $v_{st}$ decision variables. In particular, a student will be scheduled to take an exam during a time slot if and only if they are enrolled in a course group that is assigned to that time slot.
    \begin{equation} 
       M_2\cdot v_{st} \geq \sum_{g\in \mathcal{G}}b_{sg}x_{gt} \quad \mbox{ for } s \in \mathcal{S}, t \in \mathcal{T} \label{const2} 
    \end{equation}
    \begin{equation}
       v_{st} \leq \sum_{g\in \mathcal{G}}b_{sg}x_{gt} \quad \mbox{ for } s \in \mathcal{S}, t \in \mathcal{T} \label{const3} 
   \end{equation}
    
    \item A course group can only be assigned to a time slot if that time slot is available ($a_t=1$).
    \begin{equation} 
        x_{gt} \leq a_t \quad \mbox{ for } g \in \mathcal{G}, t \in \mathcal{T} \label{const4}
    \end{equation}
    
     \item At most $M_1$ students can be assigned to take a final exam during any single time slot (due to space limitations, for example).
    \begin{equation} 
       \sum_{g\in \mathcal{G}} N_g x_{gt} \leq M_1 \quad \mbox{ for } t \in \mathcal{T} \label{const5}
    \end{equation}
    
    \item A course group must be assigned to a particular time slot if required ($r_{gt}=1$).
    \begin{equation} 
      x_{gt} \geq r_{gt} \quad \mbox{ for } g \in \mathcal{G}, t \in \mathcal{T} \label{const6}
    \end{equation}
    
    \item A course group cannot be assigned to a particular time slot if that assignment is forbidden ($q_{gt}=1$).
    \begin{equation} 
      x_{gt} \leq q_{gt} \quad \mbox{ for } g \in \mathcal{G}, t \in \mathcal{T} \label{const7}
    \end{equation}
    
    \item It is not desirable for a student to be enrolled in two course groups that are assigned the same time slot (called an \emph{overlap}); this requires $z_{st}^\text{overlap}=1$ in the following constraint, incurring a penalty of $\rho_\text{overlap}$ in the objective function. Students cannot be enrolled in more than two course groups that are assigned the same time slot.
    \begin{equation} 
    \sum_{g\in \mathcal{G}}b_{sg}x_{gt} \leq 1 + z_{st}^\text{overlap}  \quad \mbox{ for } s \in \mathcal{S}, t \in \mathcal{T} \label{const8}
    \end{equation}
    
    \item It is not desirable for a student to be scheduled to take two exams in back-to-back time slots; this requires $z_{st_0}^\text{B2B}=1$ in the following constraint (where $t_0$ is the first of the two time slots), incurring a penalty of $\rho_\text{B2B}$ in the objective function. This penalty is incurred for every occurrence of back-to-back exams.
    \begin{equation} 
    v_{st_0}+v_{st_1} \leq 1 + z_{st_0}^\text{B2B}  \quad \mbox{ for } s \in \mathcal{S}, \{t_0,t_1\} \in \mathcal{T}_\text{B2B} \label{const9}
    \end{equation}
    
     \item It is not desirable for a student to be scheduled to take a night exam followed by a morning exam the next day; this requires $z_{st_0}^\text{PMtoAM}=1$ in the following constraint (where $t_0$ is the night exam time slot), incurring a penalty of $\rho_\text{PMtoAM}$ in the objective function. This penalty is incurred for every occurrence of night-to-morning exams.
    \begin{equation} 
    v_{st_0}+v_{st_1} \leq 1 + z_{st_0}^\text{PMtoAM}  \quad \mbox{ for } s \in \mathcal{S}, \{t_0,t_1\} \in \mathcal{T}_\text{PMtoAM} \label{const10}
    \end{equation}
    
    \item It is not desirable for a student to be scheduled to take three exams within a 24-hour time period; this requires $z_{s}^\text{3in24}=1$ in the following constraint, incurring a penalty of $\rho_\text{3in24}$ in the objective function. This penalty can only be incurred once per student.
    \begin{equation} 
   \sum_{t\in T} v_{st} \leq 2+z_{s}^\text{3in24}  \quad \mbox{ for } s \in \mathcal{S}, T \in \mathcal{T}_\text{3in24} \label{const11}
    \end{equation}
    
    \item It is not desirable for a student to be scheduled to take four exams within a 48-hour time period; this requires $z_{s}^\text{4in48}=1$ in the following constraint, incurring a penalty of $\rho_\text{4in48}$ in the objective function. This penalty can only be incurred once per student.
    \begin{equation} 
   \sum_{t\in T} v_{st} \leq 3 + z_{s}^\text{4in48}  \quad \mbox{ for } s \in \mathcal{S}, T \in \mathcal{T}_\text{4in48} \label{const12}
    \end{equation}
    
     \item These constraints relate the $x_{gt}$ and $w_{ft}$ decision variables. In particular a faculty member will be scheduled to give an exam during a time slot if and only if they teach a course in a course group that is assigned to that time slot.
    \begin{equation} 
       M_3\cdot w_{ft} \geq \sum_{g\in \mathcal{G}}d_{fg}x_{gt} \quad \mbox{ for } f \in \mathcal{F}, t \in \mathcal{T} \label{const13} 
    \end{equation}
    \begin{equation}
       w_{ft} \leq \sum_{g\in \mathcal{G}}d_{fg}x_{gt} \quad \mbox{ for } f \in \mathcal{F}, t \in \mathcal{T} \label{const14} 
   \end{equation}
   
    \item It is not desirable for a faculty member to teach courses in two course groups that are assigned the same time slot (called a \emph{ faculty overlap}); this requires $z_{f}^\text{facoverlap}=1$ in the following constraint, incurring a penalty of $\rho_\text{facoverlap}$ in the objective function. This penalty can only be incurred once per faculty member. Faculty members cannot teach classes in more than two course groups that are assigned the same time slot. 
    \begin{equation} 
    \sum_{g\in \mathcal{G}}d_{fg}x_{gt} \leq 1 + z_{f}^\text{facoverlap}  \quad \mbox{ for } f \in \mathcal{F}, t \in \mathcal{T} \label{const15}
    \end{equation}
    
    \item It is not desirable for a faculty member to be scheduled to give two exams in back-to-back time slots; this requires $z_{f}^\text{facB2B}=1$ in the following constraint, incurring a penalty of $\rho_\text{facB2B}$ in the objective function. This penalty can only be incurred once per faculty member.
    \begin{equation} 
    w_{ft_0}+w_{ft_1} \leq 1 + z_{f}^\text{facB2B}  \quad \mbox{ for } f \in \mathcal{F}, \{t_0,t_1\} \in \mathcal{T}_\text{B2B} \label{const16}
    \end{equation}
    
\end{itemize}

\noindent \textbf{Objective Function:} \\
We seek to assign each course group to a time slot while minimizing the penalties incurred due to undesirable aspects of the final exam schedule. 
\begin{align*}
 \sum_{s \in \mathcal{S}} \sum_{t \in \mathcal{T}}&  (\rho_\text{overlap}z_{st}^\text{overlap} +  \rho_\text{B2B}z_{st}^\text{B2B} +  \rho_\text{PMtoAM}z_{st}^\text{PMtoAM}) \\
&+\sum_{s \in \mathcal{S}} (\rho_\text{3in24}z_{s}^\text{3in24} + \rho_\text{4in48}z_{s}^\text{4in48})\\
&+\sum_{f \in \mathcal{F}} (\rho_\text{facoverlap}z_{f}^\text{facoverlap} + \rho_\text{facB2B}z_{f}^\text{facB2B})\\
\end{align*}

\noindent \textbf{Complete Formulation:} 

\noindent
\begin{flalign*}
\text{minimize}\quad &  \sum_{s \in \mathcal{S}} \sum_{t \in \mathcal{T}}  (\rho_\text{overlap}z_{st}^\text{overlap} +  \rho_\text{B2B}z_{st}^\text{B2B} +  \rho_\text{PMtoAM}z_{st}^\text{PMtoAM}) &&\\
&+\sum_{s \in \mathcal{S}} (\rho_\text{3in24}z_{s}^\text{3in24} + \rho_\text{4in48}z_{s}^\text{4in48})&&\\
&+\sum_{f \in \mathcal{F}} (\rho_\text{facoverlap}z_{f}^\text{facoverlap} + \rho_\text{facB2B}z_{f}^\text{facB2B})&&
\end{flalign*}
\begin{flalign}
\text{subject to}\quad &  \sum_{t \in \mathcal{T}} x_{gt} =1 && \mbox{ for } g \in \mathcal{G} && && && && && && \tag{\ref{const1}}\\
& M_2\cdot v_{st} \geq \sum_{g\in \mathcal{G}}b_{sg}x_{gt} && \mbox{ for } s \in \mathcal{S}, t \in \mathcal{T} && \tag{\ref{const2}}\\
& v_{st} \leq \sum_{g\in \mathcal{G}}b_{sg}x_{gt} && \mbox{ for } s \in \mathcal{S}, t \in \mathcal{T} && \tag{\ref{const3}}\\
& x_{gt} \leq a_t && \mbox{ for } g \in \mathcal{G}, t \in \mathcal{T} && \tag{\ref{const4}}\\
& \sum_{g\in \mathcal{G}} N_g x_{gt} \leq M_1 && \mbox{ for } t \in \mathcal{T}&& \tag{\ref{const5}}\\
& x_{gt} \geq r_{gt} && \mbox{ for } g \in \mathcal{G}, t \in \mathcal{T}&& \tag{\ref{const6}}\\
&  x_{gt} \leq q_{gt} && \mbox{ for } g \in \mathcal{G}, t \in \mathcal{T}&& \tag{\ref{const7}}\\
& \sum_{g\in \mathcal{G}}b_{sg}x_{gt} \leq 1 + z_{st}^\text{overlap}  && \mbox{ for } s \in \mathcal{S}, t \in \mathcal{T}&& \tag{\ref{const8}}\\ 
& v_{st_0}+v_{st_1} \leq 1 + z_{st_0}^\text{B2B}  && \mbox{ for } s \in \mathcal{S}, \{t_0,t_1\} \in \mathcal{T}_\text{B2B}&&\tag{\ref{const9}}\\
& v_{st_0}+v_{st_1} \leq 1 + z_{st_0}^\text{PMtoAM}  && \mbox{ for } s \in \mathcal{S}, \{t_0,t_1\} \in \mathcal{T}_\text{PMtoAM}&& \tag{\ref{const10}}\\
& \sum_{t\in T} v_{st} \leq 2 + z_{s}^\text{3in24}  && \mbox{ for } s \in \mathcal{S}, T \in \mathcal{T}_\text{3in24}&& \tag{\ref{const11}}\\ 
& \sum_{t\in T} v_{st} \leq 3 + z_{s}^\text{4in48}  && \mbox{ for } s \in \mathcal{S}, T \in \mathcal{T}_\text{4in48} && \tag{\ref{const12}}\\
& M_3\cdot w_{ft} \geq \sum_{g\in \mathcal{G}}d_{fg}x_{gt} && \mbox{ for } f \in \mathcal{F}, t \in \mathcal{T} && \tag{\ref{const13}}\\
& w_{ft} \leq \sum_{g\in \mathcal{G}}d_{fg}x_{gt} && \mbox{ for } f \in \mathcal{F}, t \in \mathcal{T}&& \tag{\ref{const14}}\\
& \sum_{g\in \mathcal{G}}d_{fg}x_{gt} \leq 1 + z_{f}^\text{facoverlap} && \mbox{ for } f \in \mathcal{F}, t \in \mathcal{T} && \tag{\ref{const15}}\\
& w_{ft_0}+w_{ft_1} \leq 1 + z_{f}^\text{facB2B}  && \mbox{ for } f \in \mathcal{F}, \{t_0,t_1\} \in \mathcal{T}_\text{B2B} && \tag{\ref{const16}}\\
& x_{gt}\in\{0,1\} && \mbox{ for }  g \in \mathcal{G}, t \in \mathcal{T}&& \label{const17}\\
& v_{st},z_{st}^\text{overlap},z_{st}^\text{B2B},z_{st}^\text{PMtoAM}\in\{0,1\}&& \mbox{ for }  s \in \mathcal{S}, t \in \mathcal{T}&& \label{const18}\\
& w_{ft}\in\{0,1\}&& \mbox{ for }  f \in \mathcal{F}, t \in \mathcal{T}&& \label{const19}\\
& z_{s}^\text{3in24},z_{s}^\text{4in48}\in\{0,1\} && \mbox{ for }  s \in \mathcal{S} && \label{const20}\\
& z_{f}^\text{facoverlap},z_{f}^\text{facB2B}\in\{0,1\} && \mbox{ for }  f \in \mathcal{F} \label{const21}&&
\end{flalign}

\section{Two-Phase Heuristic}\label{sec:heu}

\noindent \textbf{Two-Phase Approach:}

The two-phase optimization approach outlined in Section \ref{sec:IP} works as follows. In phase 1, we only consider a subset of course groups $\mathcal{G}'\subset \mathcal{G},$ typically consisting of the 17 to 21 course groups with the largest student enrollment, as well as those course groups that are required to be assigned to a certain time slot (i.e., all $g\in\mathcal{G}$ for which $r_{gt}=1$ for some $t\in\mathcal{T}$). We also further reduce the problem by only considering penalties for overlapping, back-to-back, and night-to-morning exams. Specifically then, the phase 1 formulation has constraints \eqref{const1}--\eqref{const10} and \eqref{const17}--\eqref{const18}, substituting $\mathcal{G}'$ for $\mathcal{G}$ as necessary, and the objective is to minimize
\begin{equation*}
     \sum_{s \in \mathcal{S}} \sum_{t \in \mathcal{T}}  (\rho_\text{overlap}z_{st}^\text{overlap} +  \rho_\text{B2B}z_{st}^\text{B2B} +  \rho_\text{PMtoAM}z_{st}^\text{PMtoAM}).
\end{equation*}
This formulation is solved with a short time limit, which is gradually increased as better incumbent solutions are found.

Phase 2 then solves the complete IP given above, while fixing the course groups from the phase 1 solution (i.e., for all $g\in\mathcal{G}',$ setting $r_{gt}=1$ whenever $x^{*}_{gt}=1,$ where $\boldsymbol{x}^{*}$ is the final solution from phase 1).  It usually runs to optimality, but we set a cap of 4 hours to ensure that the registrar can run a portfolio overnight.

\finishFEappendix

\section*{Acknowledgements}
We are especially grateful to Tim Kracker, Bucknell University Registrar, and to Vince Pellegrini, Assistant Registrar for Academic Scheduling.  Vince has a near-omniscient understanding of exam scheduling at Bucknell, and facilitated combining his domain expertise with Operations Research tools.  In addition to collaborating extensively with the Registrar's Office, we worked closely with several incredible members of Bucknell's Library \& Information Technology.  These staff members helped us  build our server and protect confidential data.  Specifically, we worked closely with and are grateful to Bucknell's Data Analytics Architect Mike Latorre, Associate Director of Enterprise Technologies Systems \& Operations Jennifer Harper, and Senior Cloud Systems Engineer Wade Hutchison.  Finally, we thank the Bucknell Provost's Office for supporting our work through an Interdisciplinary Collaborations Grant.

\bibliographystyle{abbrvnat}
\bibliography{bibliog} 

\begin{thebibliography}{34}
\providecommand{\natexlab}[1]{#1}
\providecommand{\url}[1]{\texttt{#1}}
\expandafter\ifx\csname urlstyle\endcsname\relax
  \providecommand{\doi}[1]{doi: #1}\else
  \providecommand{\doi}{doi: \begingroup \urlstyle{rm}\Url}\fi

\bibitem[Akbarzadeh et~al.(2022)Akbarzadeh, Wouters, Sys, and Maenhout]{Akb22}
B.~Akbarzadeh, J.~Wouters, C.~Sys, and B.~Maenhout.
\newblock The scheduling of medical students at {Ghent University}.
\newblock \emph{INFORMS Journal on Applied Analytics}, 52\penalty0
  (4):\penalty0 303--323, 2022.

\bibitem[Al-Hawari et~al.(2020)Al-Hawari, Al-Ashi, Abawi, and Alouneh]{Alh20}
F.~Al-Hawari, M.~Al-Ashi, F.~Abawi, and S.~Alouneh.
\newblock A practical three-phase {ILP} approach for solving the examination
  timetabling problem.
\newblock \emph{International Transactions in Operational Research},
  27\penalty0 (2):\penalty0 924--944, 2020.

\bibitem[Aldeeb et~al.(2019)Aldeeb, Al-Betar, Abdelmajeed, Younes, AlKenani,
  Alomoush, Alissa, and Alqahtani]{Ald19}
B.~A. Aldeeb, M.~A. Al-Betar, A.~O. Abdelmajeed, M.~J. Younes, M.~AlKenani,
  W.~Alomoush, K.~A. Alissa, and M.~A. Alqahtani.
\newblock A comprehensive review of uncapacitated university examination
  timetabling problem.
\newblock \emph{International Journal of Applied Engineering Research},
  14\penalty0 (24):\penalty0 4524--4547, 2019.

\bibitem[Avella et~al.(2022)Avella, Boccia, Mannino, and Viglione]{Ave22}
P.~Avella, M.~Boccia, C.~Mannino, and S.~Viglione.
\newblock Practice summary: Solving the external candidates exam schedule in
  {Norway}.
\newblock \emph{INFORMS Journal on Applied Analytics}, 52\penalty0
  (2):\penalty0 226--231, 2022.

\bibitem[Ba{\c{s}}ar and Kul(2022)]{Bac22}
M.~S. Ba{\c{s}}ar and S.~Kul.
\newblock A student-based central exam scheduling model using {A*} algorithm.
\newblock \emph{Open Computer Science}, 12\penalty0 (1):\penalty0 181--190,
  2022.

\bibitem[Broder(1964)]{Bro64}
S.~Broder.
\newblock Final examination scheduling.
\newblock \emph{Communications of the ACM}, 7\penalty0 (8):\penalty0 494--498,
  1964.

\bibitem[Burke et~al.(1996)Burke, Elliman, Ford, and Weare]{Bur96b}
E.~Burke, D.~Elliman, P.~Ford, and R.~Weare.
\newblock Examination timetabling in {B}ritish universities: A survey.
\newblock In \emph{Practice and Theory of Automated Timetabling: First
  International Conference Edinburgh, UK, August 29--September 1, 1995 Selected
  Papers 1}, pages 76--90. Springer, 1996.

\bibitem[Carter(1986)]{Car86}
M.~W. Carter.
\newblock A survey of practical applications of examination timetabling
  algorithms.
\newblock \emph{Operations Research}, 34\penalty0 (2):\penalty0 193--202, 1986.

\bibitem[Carter and Laporte(1995)]{Car95}
M.~W. Carter and G.~Laporte.
\newblock Recent developments in practical examination timetabling.
\newblock In \emph{International Conference on the Practice and Theory of
  Automated Timetabling}, pages 1--21. Springer, 1995.

\bibitem[Carter et~al.(1996)Carter, Laporte, and Lee]{Car96}
M.~W. Carter, G.~Laporte, and S.~Y. Lee.
\newblock Examination timetabling: Algorithmic strategies and applications.
\newblock \emph{Journal of the Operational Research Society}, 47\penalty0
  (3):\penalty0 373--383, 1996.

\bibitem[Ceschia et~al.(2023)Ceschia, Di~Gaspero, and Schaerf]{Ces23}
S.~Ceschia, L.~Di~Gaspero, and A.~Schaerf.
\newblock Educational timetabling: Problems, benchmarks, and state-of-the-art
  results.
\newblock \emph{European Journal of Operational Research}, 308\penalty0
  (1):\penalty0 1--18, 2023.

\bibitem[Christou et~al.(2024)Christou, Vagianou, and Vardoulias]{Chr24}
I.~T. Christou, E.~Vagianou, and G.~Vardoulias.
\newblock Planning courses for student success at the {American College of
  Greece}.
\newblock \emph{INFORMS Journal on Applied Analytics}, 54\penalty0
  (4):\penalty0 365--379, 2024.

\bibitem[Cowling et~al.(2002)Cowling, Kendall, and Hussin]{Cow02}
P.~Cowling, G.~Kendall, and N.~M. Hussin.
\newblock A survey and case study of practical examination timetabling
  problems.
\newblock In \emph{Proceedings of the 4th International Conference on the
  Practice and Theory of Automated Timetabling PATAT02}, pages 258--261, 2002.

\bibitem[Dimopoulou and Miliotis(2001)]{Dim01}
M.~Dimopoulou and P.~Miliotis.
\newblock Implementation of a university course and examination timetabling
  system.
\newblock \emph{European Journal of Operational Research}, 130\penalty0
  (1):\penalty0 202--213, 2001.

\bibitem[Garc{\'\i}a-S{\'a}nchez et~al.(2019)Garc{\'\i}a-S{\'a}nchez,
  Hern{\'a}ndez, Caro, and Jim{\'e}nez]{Gar19}
{\'A}.~Garc{\'\i}a-S{\'a}nchez, A.~Hern{\'a}ndez, E.~Caro, and G.~Jim{\'e}nez.
\newblock {Universidad Polit{\'e}cnica de Madrid} uses integer programming for
  scheduling weekly assessment activities.
\newblock \emph{INFORMS Journal on Applied Analytics}, 49\penalty0
  (2):\penalty0 104--116, 2019.

\bibitem[Garey and Johnson(1990)]{Gar79}
M.~R. Garey and D.~S. Johnson.
\newblock \emph{Computers and Intractability; A Guide to the Theory of
  NP-Completeness}.
\newblock W. H. Freeman \& Co., USA, 1990.
\newblock ISBN 0716710455.

\bibitem[Gashgari et~al.(2018)Gashgari, Alhashimi, Obaid, Palaniswamy, Aljawi,
  and Alamoudi]{Gas18}
R.~Gashgari, L.~Alhashimi, R.~Obaid, T.~Palaniswamy, L.~Aljawi, and
  A.~Alamoudi.
\newblock A survey on exam scheduling techniques.
\newblock In \emph{2018 1st International Conference on Computer Applications
  \& Information Security (ICCAIS)}, pages 1--5. IEEE, 2018.

\bibitem[Gonzalez et~al.(2018)Gonzalez, Richards, and Newman]{Gon18}
G.~Gonzalez, C.~Richards, and A.~Newman.
\newblock Optimal course scheduling for {United States Air Force Academy}
  cadets.
\newblock \emph{Interfaces}, 48\penalty0 (3):\penalty0 217--234, 2018.

\bibitem[G{\"u}ler and Gecici(2020)]{Gul20}
M.~G. G{\"u}ler and E.~Gecici.
\newblock A spreadsheet-based decision support system for examination
  timetabling.
\newblock \emph{Turkish Journal Of Electrical Engineering And Computer
  Sciences}, 28\penalty0 (3):\penalty0 1584--1598, 2020.

\bibitem[Laporte and Desroches(1984)]{Lap84}
G.~Laporte and S.~Desroches.
\newblock Examination timetabling by computer.
\newblock \emph{Computers \& Operations Research}, 11\penalty0 (4):\penalty0
  351--360, 1984.

\bibitem[Lotfi and Cerveny(1991)]{Lot91}
V.~Lotfi and R.~Cerveny.
\newblock A final-exam-scheduling package.
\newblock \emph{Journal of the Operational Research Society}, 42\penalty0
  (3):\penalty0 205--216, 1991.

\bibitem[McCollum et~al.(2007)McCollum, McMullan, Burke, Parkes, and Qu]{Mcc07}
B.~McCollum, P.~McMullan, E.~K. Burke, A.~J. Parkes, and R.~Qu.
\newblock The second international timetabling competition: Examination
  timetabling track.
\newblock Technical report, Technical Report QUB/IEEE/Tech/ITC2007/-Exam/v4.
  0/17, Queen’s University, 2007.

\bibitem[McCollum et~al.(2012)McCollum, McMullan, Parkes, Burke, and Qu]{Mcc12}
B.~McCollum, P.~McMullan, A.~J. Parkes, E.~K. Burke, and R.~Qu.
\newblock A new model for automated examination timetabling.
\newblock \emph{Annals of Operations Research}, 194:\penalty0 291--315, 2012.

\bibitem[Mehta(1981)]{Meh81}
N.~K. Mehta.
\newblock The application of a graph coloring method to an examination
  scheduling problem.
\newblock \emph{Interfaces}, 11\penalty0 (5):\penalty0 57--65, 1981.

\bibitem[Miyake et~al.(2025)Miyake, Snyder, and Tran]{github}
T.~Miyake, L.~Snyder, and V.~Tran.
\newblock Final exam schedule optimizer.
\newblock \url{https://github.com/tm032/final_exam_scheduler/tree/main}, 2025.
\newblock Accessed: 2025-09-08.

\bibitem[Muklason et~al.(2017)Muklason, Parkes, {\"O}zcan, McCollum, and
  McMullan]{Muk17}
A.~Muklason, A.~J. Parkes, E.~{\"O}zcan, B.~McCollum, and P.~McMullan.
\newblock Fairness in examination timetabling: Student preferences and extended
  formulations.
\newblock \emph{Applied Soft Computing}, 55:\penalty0 302--318, 2017.

\bibitem[Prida~Romero(1982)]{Rom82}
B.~Prida~Romero.
\newblock Examination scheduling in a large engineering school: A
  computer-assisted participative procedure.
\newblock \emph{Interfaces}, 12\penalty0 (2):\penalty0 17--24, 1982.

\bibitem[Qu et~al.(2009)Qu, Burke, McCollum, Merlot, and Lee]{Qu09}
R.~Qu, E.~K. Burke, B.~McCollum, L.~T. Merlot, and S.~Y. Lee.
\newblock A survey of search methodologies and automated system development for
  examination timetabling.
\newblock \emph{Journal of Scheduling}, 12:\penalty0 55--89, 2009.

\bibitem[Rowling(1999)]{JKR99}
J.~K. Rowling.
\newblock \emph{Harry Potter and the Prisoner of Azkaban}.
\newblock Bloomsbury, 1999.

\bibitem[Siew et~al.(2024)Siew, Sze, Goh, Kendall, Sabar, and Abdullah]{Sie24}
E.~S.~K. Siew, S.~L. Sze, S.~L. Goh, G.~Kendall, N.~R. Sabar, and S.~Abdullah.
\newblock A survey of solution methodologies for exam timetabling problems.
\newblock \emph{IEEE Access}, 12:\penalty0 41479--41498, 2024.

\bibitem[Strichman(2017)]{Str17}
O.~Strichman.
\newblock Near-optimal course scheduling at the {Technion}.
\newblock \emph{Interfaces}, 47\penalty0 (6):\penalty0 537--554, 2017.

\bibitem[Wang et~al.(2010)Wang, Bussieck, Guignard, Meeraus, and
  O’Brien]{Wan10}
S.~Wang, M.~Bussieck, M.~Guignard, A.~Meeraus, and F.~O’Brien.
\newblock Term-end exam scheduling at {United States Military Academy / West
  Point}.
\newblock \emph{Journal of Scheduling}, 13\penalty0 (4):\penalty0 375--391,
  2010.

\bibitem[Welsh and Powell(1967)]{Wel67}
D.~J. Welsh and M.~B. Powell.
\newblock An upper bound for the chromatic number of a graph and its
  application to timetabling problems.
\newblock \emph{The Computer Journal}, 10\penalty0 (1):\penalty0 85--86, 1967.

\bibitem[Ye et~al.(2024)Ye, Jovine, van Osselaer, Zhu, and Shmoys]{Ye24}
T.~Ye, A.~Jovine, W.~van Osselaer, Q.~Zhu, and D.~B. Shmoys.
\newblock {Cornell University} uses integer programming to optimize final exam
  scheduling.
\newblock \emph{arXiv preprint arXiv:2409.04959}, 2024.

\end{thebibliography}

\end{document}